\documentclass[11pt]{amsart}
\usepackage{amssymb,amsmath,amsthm,latexsym}
\usepackage{mathrsfs}
\newtheorem{lemma}{Lemma}
\newtheorem{theorem}{Theorem}

\newtheorem{proposition}[theorem]{Proposition}

\newtheorem{remark}[theorem]{Remark}
\usepackage{graphicx}
\theoremstyle{remark}
\numberwithin{equation}{section}
\usepackage{amssymb}

\begin{document}

\title{\sc Weighted Zero-sum problems over $C_3^r$}
\author{Hemar Godinho}
\address{Departamento de Matem\'{a}tica, Universidade de Bras\'{i}lia, Bras\'{i}lia-DF, Brazil}
\email{hemar@mat.unb.br}
\thanks{The first two authors were partially supported by a grant from CNPq-Brazil}
\author{AB\' ILIO LEMOS}
\address{Departamento de Matem\'{a}tica, Universidade Federal de Vi\c cosa, Vi\c cosa-MG, Brazil}
\email{abiliolemos@ufv.br}
\author{DIEGO MARQUES}
\address{Departamento de Matem\'{a}tica, Universidade de Bras\'{i}lia, Bras\'{i}lia-DF, Brazil}
\email{diego@mat.unb.br}
\thanks{The third author is partially supported by a grant from FEMAT-Brazil and CNPq-Brazil}

\begin{abstract}
Let $C_n$ be the cyclic group of order $n$ and set $s_{A}(C_n^r)$ as the
smallest integer $\ell$ such that every sequence $\mathcal{S}$ in
$C_n^r$ of length at least $\ell$ has an $A$-zero-sum subsequence of
length equal to $\exp(C_n^r)$, for $A=\{-1,1\}$. In this paper, among other 
things, we give estimates for $s_A(C_3^r)$,  and  prove that $s_A(C_{3}^{3})=9$, $s_A(C_{3}^{4})=21$ and $41\leq s_A(C_{3}^{5})\leq45$.

\end{abstract}

\maketitle

\section{Introduction}
Let $G$ be a finite abelian group (written additively), and
$\mathcal{S}$ be a finite sequence of elements of $G$ and of
length $\mathfrak{m}$. For simplicity we are going to write
$\mathcal{S}$ in a \emph{multiplicative} form
$$
\mathcal{S}=\displaystyle\prod_{i=1}^{\ell} g_{i}^{v_{i}},
$$
where $v_{i}$ represents the number of times the element $g_{i}$
appears in this sequence. Hence $\sum_{i=1}^{\ell} v_{i} =
\mathfrak{m}$.

Let $A=\{-1,1\}$. We say that a subsequence $a_{1}\cdots a_{s}$ of $\mathcal{S}$ is
an $A$-\emph{zero-sum subsequence}, if we can find $\epsilon_1,\ldots,\epsilon_s\in A$ such that
$$
\epsilon_1 a_{1} + \cdots + \epsilon_s a_{s} = 0 \;\;\; \mbox{in}
\;\;G.
$$

Here we are particularly interested in studying the behavior of $s_{A}(G)$ defined as
the smallest integer $\ell$ such that every sequence $\mathcal{S}$ of length greater than or equal to $\ell$, satisfies the condition $(s_{A})$, which states that there must exist
an $A$-zero-sum subsequence of $\mathcal{S}$ of length $\exp(G)$ (the exponent of $G$).

For this purpose, two other invariants will be defined to help us
in this study. Thus, define $\eta_{A}(G)$ as the smallest integer
 $\ell$ such that every sequence $\mathcal{S}$ of length greater than or equal to $\ell$, satisfies the
condition $(\eta_{A})$, which says that there exists an
$A$-zero-sum subsequence of $\mathcal{S}$ of length \textit{at most} $\exp(G)$. Define also $g_{A}(G)$ as the smallest integer $\ell$ such that every sequence $\mathcal{S}$ of \textit{distinct}
elements and of length greater than or equal to $\ell$, satisfies the
condition $(g_{A})$, which says that there must exist an
$A$-zero-sum subsequence of $\mathcal{S}$ of length $\exp(G)$.

The study of zero-sums is a classical area of additive number theory and goes back to the works of
Erd\"{o}s, Ginzburg and Ziv \cite{E-G-Z} and Harborth \cite{Harborth}. A very thorough survey up
to 2006 can be found on Gao-Geroldinger \cite{Gao1}, where applications of this theory are also given.

In \cite{Gry}, Grynkiewicz established a weighted version of Erd$\ddot{\mbox{o}}$s-Ginzburg-Ziv theorem, which introduced the idea of considering certain weighted subsequence sums, and Thangadurai \cite{T} presented many results on a weighted Davenport's constant and its relation to $s_A$.

For the particular weight $A=\{-1,1\}$, the best results are due to Adhikari {\it et al} \cite{Adhikari1}, where it is proved that $s_A(C_n)=n+\lfloor \log_2 n \rfloor$ (here $C_n$ is a cyclic group of order $n$) and Adhikari {\it et al} \cite{Adhikari2}, where it is proved that $s_A(C_n\times C_n )=2n-1$, when $n$ is odd. Recently, Adhikari {\it et al} proved that $s_{A}(G) = \exp(G) + |G|\log2  + O(|G|\log2 \log2 )$ when $\exp(G)$ is even and $\exp(G)\rightarrow +\infty $ (see \cite{Adhikari3}).


The aim of this paper is to give estimates for $s_A(C_n^r)$, where as usual $C_{n}^{r}=C_n\times \cdots \times C_n$ ($r$ times), and here are our results. 

\begin{theorem}\label{theorem1}
Let $A=\{-1,1\}$,  $n>1$ odd and $r\geq1$. If  $n=3$ and $r\geq 2$, or $n\geq 5$  then
$$2^{r-1}(n-1)+1\leq s_A(C_n^r)\leq (n^r-1)\left(\frac{n-1}{2}\right)+1.$$
\end{theorem}

For the case of $n=3$ we present a more detailed study
and  prove 
\begin{theorem}\label{main}
Let $A=\{-1,1\}$ and $r\geq 5$.
\begin{enumerate}
\item[(i)] If $r$ is odd then
$$s_A(C_3^r)\geq 2^{r}+ 2{r-1 \choose \frac{r-5}{2}} -1$$
\item[(ii)] If  $r$ is even, with $m=\left\lfloor \frac{3r-4}{4}\right\rfloor$, then
\begin{enumerate}
\item [(a)] If $r\equiv 2 \pmod 4$, then $s_{A}(C_3^r)\geq 2\sum_{1\leq j\leq m}{r\choose j}+ 2{r\choose \frac{r-2}{2}}+1$, where $j$ takes odd values.

\item [(b)] If $r\equiv 0\pmod 4$, then $s_{A}(C_3^r)\geq 2\sum_{1\leq j\leq m}{r\choose j}+ {r\choose \frac{r}{2}}+1$, where $j$ takes odd values.
\end{enumerate}
\end{enumerate} 
\end{theorem}

It is simple to check that $s_{A}(C_{3})=4$, and it follows from Theorem 3 in \cite{Adhikari2} that $s_{A}(C_{3}^{2})=5$. Our next result presents both exact values of $s_{A}(C_{3}^{r})$, and $r=3,4$ as well as estimates for $s_A(C_{3^a}^r)$, $r=3,4,5$, for all $a\geq1$.
\begin{theorem}\label{theorem 3}
 Let $A=\{-1,1\}$. Then
\begin{enumerate} 
\item[(i)] {$s_{A}(C_{3}^{3})=9$, $s_{A}(C_{3}^{4})=21$, $41\leq s_A(C_{3}^{5})\leq45$} 
\item[(ii)] {$s_A(C_{3^a}^3)=4\times3^a-3$, for all $a\geq1$}
\item[(iii)] {$8\times3^a-7\leq s_A(C_{3^a}^4)\leq 10\times3^a-9$, for all $a\geq1$}
\item[(iv)] {$16\times3^a-15\leq s_A(C_{3^a}^5)\leq 22\times3^a-21$, for all $a\geq1$}
\end{enumerate}
\end{theorem}


\section{Relations between the invariants $\eta_A,\ g_A$ and $s_A$}

We start by proving the following result.
\begin{lemma}\label{L1}
For $A=\left\{-1, 1\right\}$, we have 
\begin{enumerate}
 \item[(i)]$ \eta_{A}(C_{3})=2,\;\;g_{A}(C_{3})=3\;\;\mbox{and}\;\;s_{A}(C_{3})=4,$ and
\item[(ii)] $\eta_{A}(C_{3}^{r})\geq r+1$ for any $r\in\mathbb{N}$.
\end{enumerate}
\end{lemma}
\noindent \textbf{Proof.}
The proof of item (i) is very simple and will be omitted. For (ii), the proof follows from the fact that
the sequence $e_1e_2\cdots e_r$ with $e_j=(0,\ldots,1,\ldots,0)$, has no $A$-zero-sum subsequence.
\qed
\begin{proposition}\label{a3}
For $A=\left\{-1, 1\right\}$, we have
$g_{A}(C_{3}^{r})=2\eta_{A}(C_{3}^{r})-1$.
\end{proposition}
\noindent \textbf{Proof.} The case $r=1$ follows from Lemma \ref{L1}. 
Let $\mathcal{S}=\prod_{i=1}^{\mathfrak{m}}g_{i}$ of length
$\mathfrak{m}=\eta_{A}(C_{3}^{r})-1$ which does not
 satisfy the condition $(\eta_{A})$. In particular $\mathcal{S}$ has no
 $A$-zero-sum subsequences of length 1 and 2,
 that is, all elements of $\mathcal{S}$ are nonzero and distinct.
 Now, let $\mathcal{S}^{*}$ be the sequence
 $\prod_{i=1}^{\mathfrak{m}}g_{i}\prod_{i=1}^{\mathfrak{m}}(-g_{i})$. Observe that
 $\mathcal{S}^{*}$ has only distinct elements, since $\mathcal{S}$ has no $A$-zero-sum subsequences of length
 2. It is easy to see that any $A$-zero-sum of
 $\mathcal{S}^{*}$ of length $3$ is also an $A$-zero-sum of
 $\mathcal{S}$, for $A=\{-1,1\}$. Hence $g_{A}(C_{3}^{r}) \geq 2\eta_{A}(C_{3}^{r})-1$.

Let $\mathcal{S}$ be a sequence of distinct elements and of
length $\mathfrak{m}=2\eta_{A}(C_{3}^{r})-1$, and write
$$
\mathcal{S}=\displaystyle\prod_{i=1}^{t}g_{i}\displaystyle\prod_{i=1}^{t}(-g_{i})\displaystyle\prod_{i=2t+1}^{\mathfrak{m}}g_{i}
$$
 where $g_r\neq -g_s$ for $2t+1 \leq r < s \leq \mathfrak{m}$. If $t=0$, then
 $\mathcal{S}$ has no $A$-zero-sum of length 2, and  0 can
 appear at most once in $\mathcal{S}$. Let $\mathcal{S}^{*}$ be the subsequence of all nonzero 
 elements of $\mathcal{S}$, hence $|\mathcal{S}^{*}|= 2\eta_{A}(C_{3}^{r})-2 > \eta_{A}(C_{3}^{r})$, for $r\geq
 2$ (see Lemma \ref{L1}(ii)), hence it must contain an $A$-zero-sum of length 3.

For the case $t\geq 1$, we may assume $g_{j}\neq0$, for every $\jmath=2t+1,\ldots, \mathfrak{m}$ since otherwise, $g_{t}+(-g_{t})+g_{j_0}$ is $A$-zero-sum subsequence of length $3$. But now, either
 $t\geq\eta_{A}(C_{3}^{r})$,
 so that $\prod_{i=1}^{t}g_{i}$ has an $A-$zero-sum  of length 3,
 or $\mathfrak{m}-t\geq\eta_{A}(C_{3}^{r})$, so that
 $\prod_{i=1}^{t}(-g_{i})\prod_{i=2t+1}^{\mathfrak{m}}g_{i}$ has an
 $A-$zero-sum subsequence of length 3.
\qed

Here we note that by the definition of these invariants and the
proposition above, we have
\begin{equation}\label{sAgA}
s_{A}(C_{3}^{r})\geq g_{A}(C_{3}^{r})=2\eta_{A}(C_{3}^{r})-1.
\end{equation}

\begin{proposition}\label{a4}
For $A=\left\{-1, 1\right\}$, we have 
$s_{A}(C_{3}^{r})=g_{A}(C_{3}^{r})$, for $r\geq 2$.
\end{proposition}

\noindent \textbf{Proof.}  From Theorem 3 in \cite{Adhikari2} we
have $s_{A}(C_{3}^{2})=5$ and, on the other hand, the sequence
$(1, 0)(0, 1)(2, 0)(0, 2)$ does not  satisfy the condition
$(g_{A})$, hence $s_{A}(C_{3}^{2})=g_{A}(C_{3}^{2})$ (see
(\ref{sAgA})). From now on, let us consider $r\geq 3$.

Let $\mathcal{S}$ be a sequence of length
$\mathfrak{m}=s_{A}(C_{3}^{r})-1$ which does not satisfy the
condition $(s_{A})$. In particular $\mathcal{S}$ does not contain
three equal elements, since $3g=0$. If $\mathcal{S}$ contains only
distinct elements, then it does not satisfy also the condition $(g_A)$,
 and then $\mathfrak{m}\leq g_{A}(C_{3}^{r})-1$, which implies
$s_{A}(C_{3}^{r})=g_{A}(C_{3}^{r})$ (see (\ref{sAgA})). Hence, let us assume that
$\mathcal{S}$ has repeated elements and write

\begin{equation}\label{H}
\mathcal{S}= \mathcal{E}^{2}\mathcal{F} =
\prod_{i=1}^{t}g_{i}^{2}\prod_{j=2t+1}^{\mathfrak{m}}g_{j}
\end{equation}
where $g_1,\dots,g_t,g_{2t+1}, \ldots, g_{\mathfrak{m}}$ are distinct. If for
some $1\leq j\leq \mathfrak{m}$ we have $g_{j}=0$, then the
subsequence of all nonzero elements of $\mathcal{S}$ has length at
least equal to $s_{A}(C_{3}^{r})-3 \geq  2\eta_{A}(C_{3}^{r}) -4 \geq
\eta_{A}(C_{3}^{r})$ for $r\geq 3$ (see Lemma \ref{L1} (ii)). Then it must have an
$A$-zero-sum of length 2 or 3. And if the $A$-zero-sum is of
length 2, together with $g_{j}=0$ we would have an $A$-zero-sum of
length 3 in $\mathcal{S}$, contradicting the assumption that it
does not satisfy the condition $(s_{A})$.

Hence let us assume that all elements of $\mathcal{S}$ are
nonzero. Observe that we can not have $g$ in $\mathcal{E}$ and $h$
in $\mathcal{F}$ (see (\ref{H})) such that $h=-g$, for
$g+g-h=3g=0$, an $A$-zero-sum of length 3. Therefore the new sequence
$$
 \mathcal{R}=\prod_{i=1}^{t}g_{i}\prod_{i=1}^{t}(-g_{i})\prod_{i=2t+1}^{\mathfrak{m}}g_{i}
 $$
has only distinct elements, length $\mathfrak{m} =
s_{A}(C_{3}^{r})-1$, and does not satisfy the condition $(g_{A})$.
Hence $\mathfrak{m}\leq g_{A}(C_{3}^{r})-1$, and this concludes
the proof according to (\ref{sAgA}).
\qed

\section{Proof of Theorem \ref{theorem1}}
\subsection{The lower bound for $s_{A}(C_{n}^{r})$}
Let $e_1,\ldots,e_r$ be the elements of $C_n^{r}$ defined as $e_j=(0,\ldots,0,1,0,\ldots,0)$, and for every subset 
$I \subset\{1,\ldots, r\}$, of \textit{odd cardinality}, define $\mathfrak{e}_{I}=\sum_{i\in I}e_{i}$ (e.g.,
 taking $I=\{1,3,r\}$, we have $\mathfrak{e}_{I}=(1, 0, 1,
0, \dots, 0,1 )$), and let $\mathscr{I}_{m}$ be the collection of
all subsets of $\{1,\ldots,r\}$ of cardinality odd and at most
equal to $m$.

There is a natural isomorphism between the cyclic groups $C_n^{r}\cong (\mathbb{Z}/n\mathbb{Z})^{r}$, and this result here will be proved for $(\mathbb{Z}/n\mathbb{Z})^{r}$. 
 Let $\phi:\mathbb{Z} \to \mathbb{Z}/n\mathbb{Z}$ be the canonical group epimorphism, and define 
$\varphi:\mathbb{Z}^{r} \to (\mathbb{Z}/n\mathbb{Z})^{r}$ as $\varphi(a_1,\cdots,a_r)=(\phi(a_1),\cdots,\phi(a_r))$.
If $\mathcal{S}=g_1\cdots g_m$ is a sequence over the group $\mathbb{Z}^{r}$, let us denote by $\varphi(\mathcal{S})$ the sequence $\varphi(\mathcal{S})=\varphi(g_1)\cdots \varphi(g_m)$ of same length over the group $(\mathbb{Z}/n\mathbb{Z})^{r}$.

Let $e^{*}_1,\ldots,e^{*}_r$ be the canonical basis (i.e.,$e^{*}_j=(0,\ldots,0,1,0,\ldots,0)$) of the group $\mathbb{Z}^{r}$ , and define, as above
$$\mathfrak{e}^{*}_{I}=\sum_{i\in I}e_{i}^{*}$$

Now consider the sequence 
$$
\mathcal{S}=\displaystyle\prod_{I\in\mathscr{I}_r}
(\mathfrak{e}^{*}_{I})^{n-1},
$$
of length $2^{r-1}(n-1)$. We will prove that the corresponding sequence 
$$
\varphi({\mathcal{S}})=\displaystyle\prod_{I\in\mathscr{I}_r}
\mathfrak{e}_{I}^{n-1},
$$ 
has no
A-zero-sum subsequences of length $n$, which is equivalent to prove that given $A=\{-1,1\}$ and  any subsequence $\mathcal{R}=g_1\cdots g_n$ of $\mathcal{S}$, it is not possible to find 
$\epsilon_1,\ldots,\epsilon_s\in A$ such that (with an abuse of notation)
\begin{equation}\label{eq0}
\epsilon_1g_1 +\cdots + \epsilon_ng_n \equiv (0,\ldots,0)\; (\mbox{mod}\,n).
\end{equation}
Writing  $g_k=(c_1^{(k)}, \ldots,c_r^{(k)} )$, for $1\leq k \leq n$,
 it follows from (\ref{eq0}) that, for every $j\in \{1,\dots,r\}$, we have
\begin{equation}\label{eq1}
\displaystyle\sum_{k=1}^n\epsilon_{k}c_j^{(k)}\equiv 0\ (\
\mbox{mod}\ n).
\end{equation}
For every $1\leq j \leq r$, let us define the sets
$$
A_{j}=\{\ell \; | \;c_j^{(\ell)}=1\}.
$$
Since $c_j^{(\ell)}\in\{0,1\}$  and $\epsilon_j\in \{-1,1\}$ for any $j$ and any $\ell$, we must have, according to (\ref{eq1}), that either
\begin{equation}\label{Aj}
|A_{j}| = n \;\;\;\mbox{or}\;\;\; |A_{j}|\;\;\mbox{is even}.
\end{equation}
Since $g_{\ell} = \mathfrak{e}_{I_{\ell}}$, for some $I$, by the definition we have
$\sum_{j=1}^{r}c_j^{(\ell)}= |I|$  for all $\ell$, then
$$
\sum_{j=1}^{r}|A_{j}| = \sum_{j=1}^{r}\sum_{\ell=1}^{n}c_j^{(\ell)}=
 \sum_{\ell=1}^{n}\sum_{j=1}^{r}c_j^{(\ell)}=|I_1|+\cdots+|I_{n}|,
$$
an odd sum of odd numbers. Hence there exists a $\jmath_0$, such that $|A_{j_0}|=n$ (see (\ref{Aj})),
but then, it follows from (\ref{eq1}) that $\sum_{k=1}^n\epsilon_{k}c_{j_{0}}^{(k)}=n$  and therefore
$\epsilon_{1}=\cdots =\epsilon_{n}=1$. And the important consequence is that we must have $g_1=\cdots =g_{n}$,
which is impossible since in the sequence $\mathcal{S}$ no element appears more than $n-1$ times.

\begin{remark}
If we consider the sequence $\varphi({\mathcal{S}})=\prod_{I\in\mathscr{I}_r}
\mathfrak{e}_{I}$, for $n=3$, we see that this does not satisfy the
condition $(\eta_{A})$. So $\eta_{A}(C_{3}^{r})\geq 2^{r-1}+1$ for any $r\in\mathbb{N}$, which is an improvement of the item (ii) of the Lemma \ref{L1}. 
\end{remark}
 
\subsection{The upper bound for $s_{A}(C_{n}^{r})$}

Let us consider the set of elements of the group $C_{n}^{r}$ as the union $\{0\}\cup G^{+}\cup G^{-}$, where if 
$g\in G^{+}$ then $-g\in G^{-}$. And write the sequence $\mathcal{S}$ as

$$
\mathcal{S} = 0^{m}\prod_{g\in G^{+}} (g^{v_{g}(\mathcal{S})}(-g)^{v_{-g}(\mathcal{S})}).
$$
First observe that if for some $g$, $v_{g}(\mathcal{S})+v_{-g}(\mathcal{S})\geq n$, then we can find a 
 subsequence $\mathcal{R}=c_1\cdots c_{n}$ of $\mathcal{S}$, which is an $A$-zero-sum, for $A=\{-1,1\}$, and any sum of $n$ equal elements is equal to zero in $C_{n}^{r}$.
Now consider $m\geq 1$ and   $m+v_{g}(\mathcal{S})+v_{-g}(\mathcal{S}) > n,$
then we can find a subsequence $\mathcal{R}=h_1\cdots h_{t}$ of $\mathcal{S}$ of \textit{even} length $t\geq n-m$ with $h_j\in\{-g,g\}$. Since $A=\{-1,1\}$, this is an $A$-zero-sum. Hence, the subsequence $T=0^{m^{*}}\mathcal{R}$ ($m^{*}\leq m$) of $\mathcal{S}$ is an $A$-zero-sum of length $n$.

Thus assume that, for every $g$ in $\mathcal{S}$ we have $v_{g}(\mathcal{S})+v_{-g}(\mathcal{S}) \leq n-m$, which gives
$$
|\mathcal{S}| \leq \left \{\begin{array}{ccc}
m + \frac{n^{r}-1}{2}(n-m) & \mbox{if}\;\;\;m>0\;\;\; \mbox{even}\\ 
	 m-1 + \frac{n^{r}-1}{2}(n-m) & \mbox{if}\;\;\;m>0\;\;\; \mbox{odd} \\
	 \frac{n^{r}-1}{2}(n-1) & \mbox{if}\;\;\;m = 0,
\end{array}\right.
$$
for $|G^{+}|=\frac{n^{r}-1}{2}.$ We observe than in the case $m$ even $m + \frac{n^{r}-1}{2}(n-m)\leq 2 + \frac{n^{r}-1}{2}(n-2)\leq 2 + \frac{n^{r}-1}{2}(n-2)+ \frac{n^{r}-1}{2}-1$ and the equality only happens when $n=3$ and $r=1$. In any case, if $|\mathcal{S}|\geq \frac{n^{r}-1}{2}(n-1) +1$, it has a subsequence of length $n$ which is an $A$-zero-sum.
\begin{remark} 
For $n=3$, the upper bound for $s_{A}(C_{3}^{r})$ can be improved using the result of Meshulam\cite{Meshulam} as follows. According to Proposition \ref{a4},  $s_{A}(C_{3}^{r})=g_{A}(C_{3}^{r})$ for $r\geq 2$,  and it follows from the definition that $g_{A}(C_{3}^{r})\leq g(C_{3}^{r})$, where $g(C_{3}^{r})$ is the invariant $g_A(C_{3}^{r})$ with $A=\{1\}$. Now we use the Theorem 1.2 of \cite{Meshulam} to obtain $s_{A}(C_{3}^{r})=g_{A}(C_{3}^{r})\leq g(C_{3}^{r})\leq2\times 3^r/r$. 
\end{remark}

\section{Proof of Theorem 2.}

Now we turn our attention to prove the following proposition.

\begin{proposition}\label{lower}
If $r>3$ is odd and $A=\{-1,1\}$ then $\eta_{A}(C_3^r)\geq
2^{r-1}+{r-1 \choose \delta}$,
where 
\begin{equation}\label{delta}
\delta=\delta(r) =\left \{ \begin{array}{ll}
                      \frac{(r-3)}{2} & \mbox{if}\;\; r\equiv 1\pmod{4} \\
                      \frac{(r-5)}{2} & \mbox{if}\;\; r\equiv 3\pmod{4}.
                      \end{array}\right.
\end{equation}
\end{proposition}

\noindent \textbf{Proof.} We will prove this proposition by presenting an
example of a sequence of length $2^{r-1}+{r-1 \choose \delta} -1$ with no $A$-zero-sum
subsequences of length smaller or equal to 3. Let  $\ell = {r-1 \choose \delta}$, and
 consider the sequence
$$
\mathcal{S}=\mathcal{E}.\mathcal{G}=\left(\displaystyle\prod_{I\in
\mathscr{I}_{r-2}} \mathfrak{e}_{I}\right)\cdot g_1\cdots g_{\ell},
$$
with 
$$
\begin{array}{lcc}
g_1& = & (-1,\underbrace{-1,\ldots,-1}_{\delta},1,1,\ldots,1)  \\
& \vdots & \\
g_{\ell} &=& (-1,1, \ldots, 1,\underbrace{-1,\ldots,-1}_{\delta}),
\end{array}
$$
where $\mathfrak{e}_{I}$ and $\mathscr{I}_{r-2}$ are defined in the beginning of section 2. Clearly $\mathcal{S}$ has no $A$-zero-sum subsequences of length 1 or 2 and also sum or difference of two elements of $\mathcal{G}$ will never give another element of $\mathcal{G}$, for no element of $\mathcal{G}$ has zero as one of its coordinates. Now we will consider $\mathfrak{e}_{s}-\mathfrak{e}_{t}$, where $\mathfrak{e}_{s}$ and $\mathfrak{e}_{t}$ represent the $\mathfrak{e}_{I}$'s for which $s$ coordinates are equal to 1 and $t$ coordinates are equal to 1 respectively. Thus, we see that $\mathfrak{e}_{s}-\mathfrak{e}_{t}$  will never be an element of $\mathcal{G}$ since it necessarily has either a zero coordinate or it has an odd number of 1's and -1's (and $\delta +1$ is even).

 Now, if for some $s,t$ we would have 

$$
\mathfrak{e}_{s}+\mathfrak{e}_{t} = g_i,
$$
Then $\mathfrak{e}_{t},\mathfrak{e}_{s}$ would have $\delta+1$
nonzero coordinates at the same positions (to obtain $\delta+1$ coordinates -1's).
 Hence we would need to have
$$
r +(\delta+1) = s+t 
$$
Which is impossible since $s+t$ is even and $r +(\delta+1)$ is odd, for $\delta$ is odd in any of the two cases. 

Therefore, the only  possible $A$-zero-sum
subsequence of length 3 would necessarily include one element of $\mathcal{E}$
and two elements of $\mathcal{G}$.

Let  $v,w$ be elements of $\mathcal{G}$. Now
it simple to verify that (the calculations are modulo $3$)
 either $v+w$ or $v-w$ 
have two of their entries with opposite signs (for $\delta(r)<(r-1)/2$)
 and hence  either of them can not be added to an $\,\pm\mathfrak{e}_{I}$ to obtain an $A$-zero-sum, since all its nonzero entries have
  the same sign. 
\qed

\begin{proposition}\label{a9}
Let $r>4$ be even, $m=\left\lfloor \frac{3r-4}{4}\right\rfloor$ and $A=\{-1, 1\}$. Then

 $$
 \eta_{A}(C_3^r)\geq \displaystyle \sum_{{j=1}\atop{j\,odd}}^{m}{r\choose j}+ \ell(r)+1, 
 $$
where
$$
\ell(r) =\left \{ \begin{array}{cl}
                      {r \choose \frac{r-2}{2}} & \mbox{if}\;\; r\equiv 2\pmod{4} \\
                      {r \choose \frac{r}{2}}/2 & \mbox{if}\;\; r\equiv 0\pmod{4},
                      \end{array}\right.
$$
\end{proposition}
\noindent \textbf{Proof.}

Consider the sequence $\mathcal{K}=g_1\cdots g_{\tau}$ with
$$
\begin{array}{lcc}
g_1& = & (\underbrace{-1,\ldots,-1}_{\delta},1,1,\ldots,1) \\
& \vdots & \\
g_{\tau} &=& (1,1, \ldots, 1,\underbrace{-1,\ldots,-1}_{\delta})
\end{array}
$$
where
 
$$
\begin{array}{llr}
\tau =\left \{          \begin{array}{cl}
                      \ell(r) & \mbox{if}\;\; r\equiv 2\pmod{4} \\
                      2\ell(r) & \mbox{if}\;\; r\equiv 0\pmod{4},
                      \end{array}\right.

                                          & \mbox{and} &
\delta =\left \{     \begin{array}{cl}
                     \frac{r-2}{2} & \mbox{if}\;\; r\equiv 2\pmod{4} \\
                     \frac{r}{2}& \mbox{if}\;\; r\equiv 0\pmod{4},
                      \end{array}\right.
\end{array}
$$
and rearrange the elements of the sequence $\mathcal{K}$, and write it as
$$
\mathcal{K}=\prod_{i=1}^{\tau/2} g_i\,\prod_{i=1}^{\tau/2}(-g_i)=\mathcal{K}^{+}\mathcal{K}^{-}.
$$ 
It is simple to observe that if $r\equiv 2\pmod{4}$, then $\tau=\ell$ and $\mathcal{K}^{-}=\emptyset$.

Now define the sequence
$$
\mathcal{S}=\left(\displaystyle\prod_{I\in
\mathscr{I}_{m}} \mathfrak{e}_{I}\right)\mathcal{G},
$$
where $\mathcal{G}=\mathcal{K}$ if $r\equiv 2\pmod{4}$ or $\mathcal{G}=\mathcal{K}^{+}$ if $r\equiv 0\pmod{4}$, and $m=\left\lfloor \frac{3r-4}{4}\right\rfloor$, a sequence of length $|\mathcal{S}|=\displaystyle \sum_{{j=1}\atop{j\,odd}}^{m}{r\choose j}+ \ell(r)+1$.

The first important observation is that $\mathcal{S}$ has no $A$-zero-sum subsequences of length 1 or 2. 
And also sum or difference of two elements of $\mathcal{G}$ will never be another element of $\mathcal{G}$, for it necessarily will have a zero as coordinate. Also  $\mathfrak{e}_{I}-\mathfrak{e}_{J}$  will never be an element of $\mathcal{G}$ since it necessarily has either a zero coordinate or it has an odd number of 1's and -1's (and $\delta$ is even). Now, if for some $s,t$ (both defined as in the proof of the Proposition \ref{lower}) we would have 

$$
\mathfrak{e}_{s}+\mathfrak{e}_{t} =\pm g_j,\;\;\mbox{for some}\;\jmath
$$
then $\mathfrak{e}_{t},\mathfrak{e}_{s}$ would necessarily have $\delta$
nonzero coordinates at the same positions (to obtain $\delta$ coordinates -1's).
But then
$$
s+t= r+\delta \geq \frac{3r-2}{2}, \;\;\mbox{for any value of $\delta$}
$$
which is impossible since
$$
s+t \leq 2m \leq  \frac{3r-4}{2}.
$$
 Thus the only $A$-zero-sum subsequence of length 3 possible necessarily includes an element $\mathfrak{e}_{t}$ and two elements of $\mathcal{G}$.

Let $v, w$ elements of $\mathcal {G}$. First, observe that if they do not have $-1$'s in common positions, then $v+w$ has an even amount of zeros and an even amount of $-1$'s (since $r$ and $\delta$ are both even), i.e., $v+w\neq\pm\mathfrak{e}_{I}$.  If we make $v-w$ also have an even amount of nonzero coordinates, i.e., we haven't $\,\pm \mathfrak{e}_{I}$. Now, assuming that  $v,w$  have at last a $-1$ in same position, it simple to verify that (the calculations are modulo $3$) either $v+w$ or $v-w$ have two or more of their entries with opposite signs and hence either of them can not be added to an $\,\pm\mathfrak{e}_{I}$ to obtain an $A$-zero-sum, since all its nonzero entries have the same sign. 
\qed

Theorem 2 now follows from propositions \ref{a3}, \ref{a4}, \ref{lower} and \ref{a9}. 
\section{Proof of Theorem 3.}

We start by proving the following proposition.
 
\begin{proposition}\label{a5}
For $A=\left\{-1, 1\right\}$, we have
\begin{enumerate}
 \item[(i)] $\eta_{A}(C_{3}^{2})=3$;
 \item[(ii)] $\eta_{A}(C_{3}^{3})=5$;
 \item[(iii)] $\eta_{A}(C_{3}^{4})=11$.
 \item[(iv)] $21\leq\eta_{A}(C_{3}^{5})\leq23$.
 \end{enumerate}
\end{proposition}
\noindent \textbf{Proof.}
 By Propositions \ref{a3} and \ref{a4}, we have that $s_{A}(C_{3}^{r})= g_{A}(C_{3}^{r})=2\eta_{A}(C_{3}^{r})-1$, for $r>1$, and by definition, we have $g_{A}(C_{3}^{r})\leq g(C_{3}^{r})$ resulting in  $\eta_{A}(C_{3}^{r})\leq \frac{g(C_{3}^{r}) + 1}{2}$, for $r>1$. It follows from
$$
g(C_{3}^{2})=5\; \mbox{(see \cite{K})}, g(C_{3}^{3})=10, g(C_{3}^{4})=21\; \mbox{(see \cite{Kn})}, g(C_{3}^{5})=46\;\mbox{(see \cite{Edel})}, 
$$
that  $\eta_{A}(C_{3}^{2})\leq 3$, $\eta_{A}(C_{3}^{3})\leq 5$, $\eta_{A}(C_{3}^{4})\leq 11$ and $\eta_{A}(C_{3}^{5})\leq 23$. It is easy to see that the sequences $(1,0)(0,1)$ and $(1,0,0)(0,1,0)(0,0,1)(1,1,1)$ has no $A$-zero-sum of length at most three, so $\eta_{A}(C_{3}^{2})=3$ and $\eta_{A}(C_{3}^{3})=5$. It is also simple to check that following sequences of lengths 10 and 20 respectively do not satisfy the condition ($\eta_A$):
\begin{equation}\label{b}
\begin{array}{c}
(1,1,0,0)\,\cdots\,(0,0,1,1)(1,1,1,0)\,\cdots\,(0,1,1,1)\\
 \mbox{and} \\
 (1,1,0,0,0)\, \cdots \, (0,0,0,1,1)(1,1,1,0,0)\, \cdots \, (0,0,1,1,1),
\end{array}
\end{equation}
 
hence $\eta_{A}(C_{3}^{4})=11$ and $\eta_{A}(C_{3}^{5})\geq21$.

\qed

Proposition \ref{a5} together with propositions \ref{a3} and \ref{a4} gives the proof of item (i) of Theorem \ref{theorem 3}. The proof of the remaining three items is given in Proposition \ref{cor4} below.


Before going further, we need a slight modification of a result due to Gao \textit{et al} for $A=\left\{1\right\}$ in \cite{Gao}. Here we shall use it in the
case $A=\{-1,1\}$. The proof in this case is analogous to the original one, and  shall be omit it.
\begin{proposition}\label{a7}
Let $G$ be a finite abelian group, $A=\left\{-1, 1\right\}$ and $H\leq G$. Let $\mathcal{S}$ be a sequence in $G$ of length
$$
\mathfrak{m}\geq(s_{A}(H)-1)\exp(G/H) + s_{A}(G/H).
$$
Then $\mathcal{S}$ has an $A$-zero-sum subsequence
of length $\exp(H)\exp(G/H)$. In particular, if $\exp(G)=\exp(H)\exp(G/H)$, then
$$
s_{A}(G)\leq(s_{A}(H)-1)\exp(G/H) + s_{A}(G/H).
$$
\end{proposition}

\begin{proposition}\label{cor4}
For $A=\{-1,1\}$, we have
\begin{enumerate}

\item[(i)] {$s_A(C_{3^a}^3)=4\times3^a-3$, for all $a\geq1$}

\item[(ii)] {$8\times3^a-7\leq s_A(C_{3^a}^4)\leq 10\times3^a-9$, for all $a\geq1$}

\item[(iii)] {$16\times3^a-15\leq s_A(C_{3^a}^5)\leq 22\times3^a-21$, for all $a\geq1$}

\end{enumerate}
\end{proposition}
\noindent \textbf{Proof.}
It follows of (i) from Theorem \ref{theorem 3} that $s_A(C_{3}^{3})=4\times3 - 3 = 9$. Now assume
 that $s_{A}(C_{3^{a-1}}^{3})=4\cdot3^{a-1}-3$. Thus, Proposition \ref{a7} yields
  $$
  \begin{array}{lcl}
   s_{A}(C_{3^{a}}^{3})& \leq & 3\times(s_A(C_{3^{a-1}}^{3})-1) + s_A(C_{3}^{3}) \\
                       &  \leq  &  4\times3^{a}-3
   \end{array}
   $$
 On the other hand, Theorem \ref{theorem1} gives $s_{A}(C_{3^{a}}^{3})\geq 4\times3^{a}-3$, concluding the proof of (i).
 
Again by (i) from Theorem \ref{theorem 3}, we have that $s_A(C_{3}^{4})=10\times3 - 9 = 21$. Now, assume
 that $s_{A}(C_{3^{a-1}}^{4})\leq10\cdot3^{a-1}-9$. It follows from Proposition
  \ref{a7} that
  $$
  \begin{array}{lcl}
   s_{A}(C_{3^{a}}^{4})& \leq & 3\times(s_A(C_{3^{a-1}}^{4})-1) + s_A(C_{3}^{4}) \\
                       &  \leq  &  10\times3^{a}-9
   \end{array}
   $$

On the other hand, Theorem \ref{theorem1} gives the lower bound $s_{A}(C_{3^{a}}^{4})\geq 8\times3^{a}-7$, concluding the proof of (ii). The proof of item (iii) is analogous to the proof of item (ii), again  using (i) of the Theorem \ref{theorem 3} and Theorem \ref{theorem1}.
\qed

\section*{Acknowledgement}
The authors were partially supported by a grant from CNPq-Brazil. The third author is also grateful to FEMAT-Brazil for financial support. 



\begin{thebibliography}{9999}

\bibitem{Adhikari1} S. D. Adhikari, Y. G. Chen, J. B. Friedlander, S. V. Konyagin, F. Pappalardi. Contributions to zero-sum problems. Discrete Math., 306:1-10, 2006.

\bibitem{Adhikari2} S. D. Adhikari, R. Balasubramanian, F. Pappalardi, P. Rath. Some zero-sum constants with weights. Proc. Indian Acad. Sci. (Math. Sci.), 128 (2):183-188, 2008.

\bibitem{Adhikari3} S. D. Adhikari, D. J. Grynkiewicz, Zhi-Wei Sun. On weighted zero-sum sequences. \textit{arXiv:1003.2186v1 [math.CO]} 10 Mar 2010.


\bibitem{Gao} R. Chi, S. Ding, W. Gao, A. Geroldinger, W. A. Schmid. On zero-sum subsequence of restricted size. IV. Acta Math. Hungar., 107(4):337-344, 2005.

\bibitem{Edel} Y. Edel, S. Ferret, I. Landjev, L. Storme. The classification of the largest caps in $AG(5, 3)$. J. Comb. Theory, 99:95-110, 2002.

\bibitem{E-G-Z} P. Erd\"os , A. Ginzburg and A. Ziv. Theorem in the additive number theory. Bulletim Research Council Israel 10F, 41-43, 1961.

\bibitem{Gao1} W. Gao, A. Geroldinger. Zero-sum problem in finite abelian groups: A survey. Expo. Math., 24(6): 337-369, 2006.

\bibitem{Gry} D. J. Grynkiewicz. A weighted Erd$\ddot{\mbox{o}}$s-Ginzburg-Ziv theorem. Combinatorica 26,  no. 4, 445--453, 2006.

\bibitem{Harborth} H. Harborth. Ein Extremal Problem f$\ddot{\mbox{u}}$r Gitterpunkte. J. Reine Angew. Math., 262: 356-360, 1973.

\bibitem{K} A. Kemnitz. On a lattice point problem. Ars Combinatoria, 16: 151-160, 1983.

\bibitem{Kn} D. E. Knuth, \textit{Computerprogramme}, http://www-cs-faculty.stanford.edu/$\sim$knuth/programs/setset-all.w.

\bibitem{Meshulam} R. Meshulam. On subsets of finite abelian groups with no 3-term arithmetic
progressions. J. Comb. Theory, Ser. A, 71: 168-172, 1995.

\bibitem{T} R. Thangadurai. A variant of Davenport's constant. Proc. Indian Acad. Sci. (Math. Sci.), 117: 147-158, 2007.





\end{thebibliography}
\end{document}